\documentstyle{article}
\begin{document}
\title{\Large{\bf{Lie symmetry analysis of nonlinear evolution equation for description nonlinear waves in a viscoelastic tube }}}
\author{Mehdi Nadjafikhah\thanks{Corresponding author: Department of Mathematics,
Islamic Azad University, Karaj Branch, Karaj, Iran. e-mail:
m\_nadjafikhah@iust.ac.ir} \and Vahid Shirvani-Sh.\thanks{e-mail:
v.shirvani@kiau.ac.ir } }
\date{}
\maketitle
\begin{abstract}
In this paper, the Lie symmetry method is performed for the
nonlinear evolution equation for description nonlinear waves in a
viscoelastic tube. we will find one and two-dimensional optimal
system of Lie subalgebras. Furthermore, preliminary
classification of its group-invariant solutions are investigated.
\end{abstract}


\input{amssym}
\def\be{\begin{eqnarray}}
\def\ee{\end{eqnarray}}
\def\di{\displaystyle}
\def\rank{{\bf rank}}
\section{Introduction}
The theory of Lie symmetry groups of differential equations was
developed by Sophus Lie \cite{[2]}. Such Lie groups are invertible
point transformations of both the dependent and independent
variables of the differential equations. The symmetry group
methods provide an ultimate arsenal for analysis of differential
equations and is of great importance to understand and to
construct solutions of differential equations. Several
applications of Lie groups in the theory of differential equations
were discussed in the literature, the most important ones are:
reduction of order of ordinary differential equations,
construction of invariant solutions, mapping solutions to other
solutions and the detection of linearizing transformations (for
many other applications of Lie symmetries see \cite{[3]},
\cite{[4]} and \cite{[7]}).

In the present paper, we study the following fifth-order nonlinear
evolution equation
\be
u_t+auu_x+bu_{x^3}+cu_{x^4}+du_{x^5}=eu_{x^2}.\label{eq:1.0}\ee
where $a$,$b$,$c$,$d$ and $e$ are positive constants. This
equation was introduced recently by Kudryashov and Sinelshchikov
\cite{[8]} which is the generalization of the famous Kawahara
equation. By using the reductive perturbation method, they
obtained the equation (\ref{eq:1.0}). Study of nonlinear wave
processes in viscoelastic tube is the important problem such
tubes similar to large arteries (see \cite{[9]},\cite{[10]} and
\cite {[11]}).

In this paper, by using the lie point symmetry method, we will
investigate the equation (\ref{eq:1.0}) and looking the
representation of the obtained symmetry group on its Lie algebra.
We will find the preliminary classification of group-invariant
solutions and then we can reduce the equation (\ref{eq:1.0}) to an
ordinary differential equation.

This work is organized as follows. In section 2 we recall some
results needed to construct Lie point symmetries of a given system
of differential equations. In section 3, we give the general form
of a infinitesimal generator admitted by equation (\ref{eq:1.0})
and find transformed solutions. Section 4, is devoted to the
construction of the group-invariant solutions and its
classification which provides in each case reduced forms of
equation (\ref{eq:1.0}).
\section{Method of Lie Symmetries}
In this section, we recall the general procedure for determining
symmetries for any system of partial differential equations (see
\cite{[3]} and \cite{[4]}). To begin, let us consider the general
case of a nonlinear system of partial differential equations of
order $n$th in $p$ independent and $q$ dependent variables is
given as a system of equations
\be \Delta_\nu(x,u^{(n)})=0,\;\;\;\;\; \nu=1,\cdots,l,
\label{eq:2.1} \ee
involving $x = (x^1,\cdots, x^p)$, $u = (u^1,\cdots,u^q)$ and the
derivatives of $u$ with respect to $x$ up to $n$, where $u^{(n)}$
represents all the derivatives of $u$ of all orders from $0$ to
$n$. We consider a one-parameter Lie group of infinitesimal
transformations acting on the independent and dependent variables
of the system (\ref{eq:2.1})
\be \tilde{x}^i &=& x^i+s \xi^i(x,u)+O(s^2), \hspace{1cm}
i=1\cdots,p,\nonumber \\[-2mm] \label{eq:2.2}\\[-2mm] \tilde{u}^j &=& u^j+s
\eta^j(x,u)+O(s^2), \hspace{9mm} j=1\cdots,q, \nonumber
 \ee
where $s$ is the parameter of the transformation and $\xi^i$,
$\eta^j$ are the infinitesimals of the transformations for the
independent and dependent variables, respectively. The
infinitesimal generator ${\mathbf v}$ associated with the above
group of transformations can be written as
\be  {\mathbf v} = \sum_{i=1}^p\xi^i(x,u)\partial_{x^i} +
\sum_{j=1}^q\eta^j(x,u)\partial_{u^j}. \label{eq:2.4} \ee
A symmetry of a differential equation is a transformation which
maps solutions of the equation to other solutions. The invariance
of the system (\ref{eq:2.1}) under the infinitesimal
transformations leads to the invariance conditions (Theorem 2.36
of \cite{[3]})
\be \textrm{Pr}^{(n)}{\mathbf
v}\big[\Delta_\nu(x,u^{(n)})\big]=0,\;\;\;\;\;
\nu=1,\cdots,l,\;\;\;\;\mbox{whenever}\;\;\;\;\;\Delta_\nu(x,u^{(n)})
=0, \label{eq:2.5} \ee
where $\textrm{Pr}^{(n)}$ is called the $n^{th}$ order
prolongation of the infinitesimal generator given by
\be \textrm{Pr}^{(n)}{\mathbf v}= {\mathbf
v}+\sum^q_{\alpha=1}\sum_J
\varphi^J_\alpha(x,u^{(n)})\partial_{u^\alpha_J},\label{eq:2.6}
\ee
where $J=(j_1,\cdots,j_k)$, $1\leq j_k\leq p$, $1\leq k\leq n$ and
the sum is over all $J$'s of order $0<\# J\leq n$. If $\#J=k$, the
coefficient $\varphi_J^\alpha$ of $\partial_{u_J^\alpha}$ will
only depend on $k$-th and lower order derivatives of $u$, and
\be \varphi_\alpha^J(x,u^{(n)})=D_J(\varphi_\alpha - \sum_{i=1}^p
\xi^iu_i^\alpha) + \sum_{i=1}^p\xi^iu^\alpha_{J,i}, \label{eq:2.7}
\ee
where $u_i^\alpha:=\partial u^\alpha/\partial x^i$ and
$u_{J,i}^\alpha := \partial u_J^\alpha/\partial x^i$.

\medskip One of the most important properties of these
infinitesimal symmetries is that they form a Lie algebra under the
usual Lie bracket.
\section{Lie symmetries of the equation (\ref{eq:1.0}) }
We consider the one parameter Lie group of infinitesimal
transformations on $(x^1=x,x^2=t,u^1=u)$,
\be \tilde{x} &=& x+s\xi(x,t,u)+O(s^2),\nonumber\\
\tilde{t} &=& x+s\eta(x,t,u)+O(s^2),\label{eq:3.1}\\
\tilde{u} &=& x+s\varphi(x,t,u)+O(s^2),\nonumber \ee
where $s$ is the group parameter and $\xi^1=\xi$, $\xi^2=\eta$ and
$\varphi^1=\varphi$ are the infinitesimals of the transformations
for the independent and dependent variables, respectively. The
associated vector field is of the form:
\be {\mathbf
v}=\xi(x,t,u)\partial_x+\eta(x,t,u)\partial_t+\varphi(x,t,u)\partial_u.
\label{eq:3.2}\ee
and, by (\ref{eq:2.6}) its third prolongation is
\be \textrm{Pr}^{(5)}{\mathbf v} &=& {\mathbf v}+
\varphi^x\,\partial_{u_x}+\varphi^t\,\partial_{u_t}+\varphi^{x^2}\,\partial_{u_{x^2}}
+\varphi^{xt}\,\partial_{u_{xt}}+\cdots \nonumber\\
&&
\cdots+\varphi^{t^2}\,\partial_{u_{t^2}}+\varphi^{xt^4}\,\partial_{u_{xt^4}}+\varphi^{t^5}\,\partial_{u_{t^5}}.
\label{eq:3.2-1}\ee
where, for instance by (\ref{eq:2.7}) we have
\be
\varphi^x&=&D_x(\varphi-\xi\,u_x-\eta\,u_t)+\xi\,u_{x^2}+\eta\,u_{xt},\nonumber\\
\varphi^t&=&D_t(\varphi-\xi\,u_x-\eta\,u_{t})+\xi\,u_{xt}+\eta\,u_{t^2},\nonumber\\
&& \ldots \label{eq:3.3} \\
\varphi^{t^5}&=&D^5_{x}(\varphi-\xi\,u_x-\eta\,u_t)+\xi\,u_{x^5t}+\eta\,u_{t^5},\nonumber
\ee
where $D_x$ and $D_t$ are the total derivatives with respect to
$x$ and $t$ respectively.
By (\ref{eq:2.5}) the vector field ${\mathbf v}$ generates a one
parameter symmetry group of equation (\ref{eq:1.0}) if and only if
\be \left. \begin{array}{l} \di \textrm{Pr}^{(5)}{\mathbf
v}[u_t+auu_x+bu_{x^3}+cu_{x^4}+du_{x^5}-eu_{x^2}]=0,\\[5mm]
\di \mbox{whenever} \hspace{.5cm}
u_t+auu_x+bu_{x^3}+cu_{x^4}+du_{x^5}-eu_{x^2}=0.
\end{array}\right. \label{eq:3.3-1} \ee
The condition (\ref{eq:3.3-1}) is equivalent to
\be \left. \begin{array}{l} \di au_{x}\varphi+
au\varphi^x+\varphi^t-e\varphi^{x^2}+b\varphi^{x^3}+c\varphi^{x^4}+d\varphi^{x^5}=0,\\[5mm]
\di \mbox{whenever} \hspace{.5cm}
u_t+auu_x+bu_{x^3}+cu_{x^4}+du_{x^5}-eu_{x^2}=0.
\end{array} \right. \label{eq:3.4} \ee

Substituting (\ref{eq:3.3}) into (\ref{eq:3.4}), and equating the
coefficients of the various monomials in partial derivatives with
respect to $x$ and various power of $u$, we can find the
determining equations for the symmetry group of the equation
(\ref{eq:1.0}). Solving this equation, we get the following forms
of the coefficient functions
\be
 \xi=c_2at+c_3, \quad \eta=c_1, \quad  \varphi=c_2. \label{eq:3.5} \ee
where $c_1$, $c_2$ and $c_3$ are arbitrary constant. Thus, the Lie
algebra of infinitesimal symmetry of the equation (\ref{eq:1.0})
is spanned bye the three vector fields:
\be \textbf{v}_1=\partial_x,\quad \textbf{v}_2=\partial_t,\quad
\textbf{v}_3=t\,\partial_x+\frac{1}{a}\partial_u.
\label{eq:3.6}\ee
The commutation relations between these vector fields are given in
the Table 1.
\begin{table}[h] \label{Tab:1}
  \caption{The commutator table} \begin{center}
\begin{tabular}{|c|c|c|c|}
  \hline
  $[{\mathbf v}_{i},{\mathbf v}_{j}]$ & ${\mathbf v}_1$ & ${\mathbf v}_2$ & ${\mathbf v}_3$ \\ \hline
  ${\mathbf v}_1$ & 0 & 0 & 0 \\ \hline
  ${\mathbf v}_2$ & 0 & 0 & ${\mathbf v}_1$ \\ \hline
  ${\mathbf v}_3$ & 0 & $-{\mathbf v}_1$ & 0 \\ \hline
\end{tabular} \end{center}
\end{table}
\paragraph{Theorem 3.1}
{ \it The Lie algebra $\pounds_{3}$ spanned by
$v_{1},v_{2},v_{3}$ is second Bianchi class type and it's
solvable and Nilpotent. }\cite{[12]}

\medskip To obtain the group transformation which is generated by
the infinitesimal generators $\textbf{v}_i$ for $i=1,2,3$ we need
to solve the three systems of first order ordinary differential
equations
\be \di \frac{d\tilde{x}(s)}{ds} &=&
\xi_i(\tilde{x}(s),\tilde{t}(s),\tilde{u}(s)), \quad
\tilde{x}(0)=x, \nonumber\\
\di \frac{d\tilde{t}(s)}{ds} &=&
\eta_i(\tilde{x}(s),\tilde{t}(s),\tilde{u}(s)), \quad
\tilde{t}(0)=t, \qquad i=1,2,3 \label{eq:3.7}\\
\di \frac{d\tilde{u}(s)}{ds} &=&
\varphi_i(\tilde{x}(s),\tilde{t}(s),\tilde{u}(s)), \quad
\tilde{u}(0)=u. \nonumber
 \ee
Exponentiating the infinitesimal symmetries of (\ref{eq:1.0}), we
get the one-parameter groups $G_i(s)$ generated by $\textbf{v}_i$
for $i=1,2,3$
\be
G_1:(t,x,u) & \longmapsto & (x+s,t,u),\nonumber\\
G_2:(t,x,u) & \longmapsto & (x,t+s,u),\label{eq:3.8} \\
G_3:(t,x,u) & \longmapsto & (x+ts,t,u+s/a).\nonumber  \ee
Consequently,
\paragraph{Theorem 3.2}
{\it If $u=f(x,t)$ is a solution of (\ref{eq:1.0}), so are the
functions
\be G_1(s)\cdot f(x,t)&=&f(x-s,t),\nonumber\\
G_2(s)\cdot f(x,t)&=&f(x,t-s),\label{eq:3.9} \\
G_3(s)\cdot f(x,t)&=&f(x-ts,t)+s/a.\nonumber \ee}
\section{Optimal system and invariant solution of (\ref{eq:1.0}) }
In this section, we obtain the optimal system and reduced forms of
the equation (\ref{eq:1.0}) by using symmetry group properties
obtained in previous section. Since the original partial
differential equation has two independent variables, then this
partial differential equation transform into the ordinary
differential equation after reduction.
\paragraph{Definition 4.1}
Let $G$ be a Lie group with Lie algebra $\goth g$. An optimal
system of $s-$parameter subgroups is a list of conjugacy
inequivalent $s-$parameter subalgebras with the property that any
other subgroup is conjugate to precisely one subgroup in the
list. Similarly, a list of $s-$parameter subalgebras forms an
optimal system if every $s-$parameter subalgebra of $\goth g$ is
equivalent to a unique member of the list under some element of
the adjoint representation: $\overline{\goth h}={\mathrm
Ad}(g({\goth h}))$.\cite{[3]}
\paragraph{Theorem 4.2}
{\it Let $H$ and $\overline{H}$ be connected s-dimensional Lie
subgroups of the Lie group $G$ with corresponding Lie subalgebras
$\goth h$ and $\overline{\goth h}$ of the Lie algebra $\goth g$ of
$G$. Then $\overline{H}$=$gHg^{-1}$ are conjugate subgroups if and
only if $\overline{\goth h}={\mathrm Ad}(g({\goth h}))$ are
conjugate subalgebras. }\cite{[3]}

\medskip By theorem (4.2), the problem of finding an optimal system of
subgroups is equivalent to that of finding an optimal system of
subalgebras. For one-dimensional subalgebras, this classification
problem is essentially the same as the problem of classifying the
orbits of the adjoint representation, since each one-dimensional
subalgebra is determined by nonzero vector in $\goth g$. This
problem is attacked by the na\"{\i}ve approach of taking a general
element ${\mathbf V}$ in $\goth g$ and subjecting it to various
adjoint transformation so as to "simplify" it as much as
possible. Thus we will deal with th construction of the optimal
system of subalgebras of $\goth g$.

\medskip To compute the adjoint representation, we use the Lie
series
\be {\mathrm Ad}(\exp(\varepsilon{\mathbf v}_i){\mathbf v}_j) =
{\mathbf v}_j-\varepsilon[{\mathbf v}_i,{\mathbf
v}_j]+\frac{\varepsilon^2}{2}[{\mathbf v}_i,[{\mathbf
v}_i,{\mathbf v}_j]]-\cdots,\label{eq:4.1} \ee
where $[{\mathbf v}_i,{\mathbf v}_j]$ is the commutator for the
Lie algebra, $\varepsilon$ is a parameter, and $i,j=1,2,3$. Then
we have the Table 2.
\begin{table}[h] \label{Tab:2}
\caption{Adjoint representation table of the infinitesimal
generators ${\mathbf v}_i$}
\begin{center}
\begin{tabular}{|c|c|c|c|}
  \hline
  ${\mathrm Ad}(\exp(\varepsilon {\mathbf v}_i)) {\mathbf v}_j$ & ${\mathbf v}_1$ & ${\mathbf v}_2$ & ${\mathbf v}_3$ \\\hline
  ${\mathbf v}_1$ & ${\mathbf v}_1$ & ${\mathbf v}_2$ & ${\mathbf v}_3$ \\ \hline
  ${\mathbf v}_2$ & ${\mathbf v}_1$ & ${\mathbf v}_2$ & ${\mathbf v}_3-\varepsilon{\mathbf v}_1$ \\ \hline
  ${\mathbf v}_3$ & ${\mathbf v}_1$ & ${\mathbf v}_2+\varepsilon{\mathbf v}_1$& ${\mathbf v}_3$ \\ \hline
\end{tabular} \end{center}
\end{table}

\paragraph{Theorem 4.3}
{ \it An optimal system of one-dimensional Lie algebras of the
equation (\ref{eq:1.0}) is provided by }
1) $\;\textbf{v}_2$, \quad\quad 2)
$\;\textbf{v}_3+\alpha\textbf{v}_2$ \hfill\ \mbox{ }
\paragraph{Proof:}
Consider the symmetry algebra $\goth g$ of the equation
(\ref{eq:1.0}) whose adjoint representation was determined in
table 2 and
\be {\mathbf V}=a_{1}{\mathbf v}_1+a_{2}{\mathbf v}_{2}+
a_{3}{\mathbf v}_{3}.\label{eq:4.2} \ee
is a nonzero vector field in $\goth g$. We will simplify as many
of the coefficients $a_{i}$ as possible through judicious
applications of adjoint maps to ${\mathbf V}$.
Suppose first that $a_{3}\neq0$. Scaling ${\mathbf V}$ if
necessary we can assume that $a_{3}=1$. Referring to table 2, if
we act on such a ${\mathbf V}$ by $Ad(\exp(a_{1}{\mathbf v}_2))$,
we can make the coefficient of ${\mathbf v}_1$ vanish and the
vector field ${\mathbf V}$ takes the form
\be {\mathbf V'}=Ad(\exp(a_{1}{\mathbf v}_2)){\mathbf
V}=a'_{2}{\mathbf v}_2+{\mathbf v}_3.\label{eq:4.3} \ee
for certain scalar ${a'}_2$. So, depending on the sign of
${a'}_2$, we can make the coefficient of ${\mathbf v}_2$ either
+1, -1 or 0. In other words, every one-dimensional subalgebra
generated by a ${\mathbf V}$ with $a_{3}\neq0$ is equivalent to
one spanned by either  $\textbf{v}_3+\textbf{v}_2$,
$\textbf{v}_3-\textbf{v}_2$ or $\textbf{v}_3$.

\medskip The remaining one-dimensional subalgebras are spanned by
vectors of the above form with $a_{3}=0$. If $a_{2}\neq0$, we
scale to make $a_{2}=1$, and then the vector field ${\mathbf V}$
takes the form
\be {\mathbf V''}={a''}_{1}{\mathbf v}_1+{\mathbf
v}_2.\label{eq:4.4} \ee
for certain scalar ${a''}_1$. Similarly we can vanish ${a''}_1$,
so every one-dimensional subalgebra generated by a ${\mathbf V}$
with $a_{3}=0$ is equivalent to the subalgebra spanned by
$\textbf{v}_2$.\hspace{11cm} $\square$
\paragraph{Theorem 4.4 }
{ \it An optimal system of two-dimensional Lie algebras of the
equation (\ref{eq:1.0}) is provided by }

\hspace{3cm}
$<\alpha\textbf{v}_2+\textbf{v}_3,\beta\textbf{v}_1+\gamma\textbf{v}_3>$\nonumber\\

Symmetry group method will be applied to the (\ref{eq:1.0}) to be
connected directly to some order differential equations. To do
this, a particular linear combinations of infinitesimals are
considered and their corresponding invariants are determined.

The equation (\ref{eq:1.0}) is expressed in the coordinates
$(x,t,u)$, so to reduce this equation is to search for its form
in specific coordinates. Those coordinates will be constructed by
searching for independent invariants $(\chi,\zeta)$ corresponding
to the infinitesimal generator. So using the chain rule, the
expression of the equation in the new coordinate allows us to the
reduced equation.

In what follows, we begin the reduction process of equation
(\ref{eq:1.0}).

\paragraph{4.5 Galilean-Invariant Solutions. }
First, consider
$\textbf{v}_3=t\,\partial_x+\frac{1}{a}\partial_u$. To determine
independent invariants $I$, we need to solve the first partial
differential equations $\textbf{v}_i$(I)=0, that is
invariants $\zeta$ and $\chi$ can be found by integrating the
corresponding characteristic system, which is
\be \frac{dt}{0}=\frac{dx}{t}=\frac{a\,du}{1}. \label{eq:4.3-1}
\ee
The obtained solution are given by
\be \chi=t,\qquad \zeta=u-\frac{x}{a\,t}. \label{eq:4.4-1} \ee
Therefore, a solution of our equation in this case is
\be u=f(x,\chi,\zeta)=\zeta+\frac{x}{a\,t}.\label{eq:4.5}\ee
The derivatives of $u$ are given in terms of $\zeta$ and $\chi$ as
\be u_x=\frac{1}{a\,t},\quad
u_{x^2}=u_{x^3}=u_{x^4}=u_{x^5}=0,\quad
u_t=\zeta_{\chi}-\frac{1}{a\,t^2}\,x.\label{eq:4.6}\ee
Substituting (\ref{eq:4.6}) into the equation (\ref{eq:1.0}), we
obtain the order ordinary differential equation
\be \zeta_{\chi}+\frac{1}{\chi}\,\zeta=0.\label{eq:4.7} \ee
The solution of this equation is $\zeta=\frac{c_1}{\chi}$.
Consequently, we obtain that
\be u(x,t)=\frac{x+a\,c_1}{a\,t}.\label{eq:4.8} \ee

\paragraph{4.6 Travelling wave solutions.}
The invariants of
$\textbf{v}_2+c_0\,\textbf{v}_1=c_0\,\partial_x+\partial_t$ are
$\chi=x-c_0\,t$ and $\zeta=u$ so the reduced form of equation
(\ref{eq:1.0}) is
\be
-c_0\,\zeta_{\chi}+a\,\zeta\,\zeta_{\chi}+b\,\zeta_{\chi^3}+c\,\zeta_{\chi^4}+d\,\zeta_{\chi^5}-e\,\zeta_{\chi^2}=0.\label{eq:4.9}
\ee
The family of the periodic solution for Eq.(\ref{eq:4.9}) when
$a=1$ takes the following form (see \cite{[8]})
\be
\zeta=a_0+A\,sn^4\{m\,\chi,k\}+B\,sn\{m\,\chi,k\}\,\frac{d}{d\chi}sn\{m\,\chi,k\}.\label{eq:4.10}\ee
where $sn\{m\,\chi,k\}$ is Jacobi elliptic function.
\paragraph{4.7}
The invariants of
$\textbf{v}_3+\beta\textbf{v}_2=t\,\partial_x+\beta\partial_t+\frac{1}{a}\partial_u$
 are $\chi=x-\frac{t^2}{2\beta}$ and $\zeta=u-\frac{t}{a\beta}$ so the reduced form
of equation (\ref{eq:1.0}) is
\be
\frac{1}{a\beta}-\frac{t}{\beta}\zeta_{\chi}+a\,\zeta\,\zeta_{\chi}+b\,\zeta_{\chi^3}+c\,\zeta_{\chi^4}+d\,\zeta_{\chi^5}-e\,\zeta_{\chi^2}=0.\label{eq:4.11}
\ee
\paragraph{4.8}
The invariants of $\textbf{v}_2=\partial_t$ are $\chi=x$ and
$\zeta=u$ then the reduced form of equation (\ref{eq:1.0}) is
\be
a\,\zeta\,\zeta_{\chi}+b\,\zeta_{\chi^3}+c\,\zeta_{\chi^4}+d\,\zeta_{\chi^5}-e\,\zeta_{\chi^2}=0.\label{eq:4.12}
\ee
\paragraph{4.9}
The invariants of $\textbf{v}_1=\partial_x$ are $\chi=t$ and
$\zeta=u$ then the reduced form of equation (\ref{eq:1.0}) is
$\zeta_{\chi}=0$, then the solution of this equation is
$u(x,t)=cte$.
\section*{Acknowledgment}
This research was supported by Islamic Azad University of Karaj
bracnch.

\end{document}